\documentclass[11pt,a4paper]{article}

\usepackage{latexsym}
\usepackage{amsmath,amsthm,amssymb}
\usepackage{rotating}
\usepackage{color}
\usepackage{geometry}
\newcommand{\ex}{{\rm exp}}

\def\beq{\begin{eqnarray}}
\def\eeq{\end{eqnarray}}

\def\beqn{\begin{eqnarray*}}  
\def\eeqn{\end{eqnarray*}}

\begin{document}

\begin{center}
  {\bf On Expectation Propagation and the Probabilistic Editor \\
    in some simple mixture problems}
\end{center}

\begin{center}
{\bf Nils Lid Hjort$^1$ and Donald Michael Titterington$^2$} 
\end{center}

\centerline{\bf $^1$Department of Mathematics, University of Oslo}
\centerline{\bf $^2$School of Mathematics \& Statistics,
  University of Glasgow$^*$\footnotetext[1]{Mike Titterington passed
    away in 2023, at the age of 77; this is the October 2010
    version of a paper we collaborated on then and (still)
    planned to extend before submitting to a journal}}

\medskip
\centerline{\bf October 2010}

\renewcommand{\baselinestretch}{1.5}\small\normalsize

\begin{abstract}
  \vspace*{0.5cm}
As for other latent-variable problems, exact Bayesian analysis is typically not
practicable for mixture problems and approximate methods have been developed. 
Variational Bayes tends to produce approximate posterior distributions for parameters
that are too tightly concentrated in having variances that are too small. The paper
identifies a few mixture problems in which Expectation Propagation and
variations thereof lead to approximate posterior distributions that
asymptotically exhibit `correct' variances and therefore stand to provide reliable interval
estimates for the unknown parameter or parameters.

  \noindent {\em Some key words}: Assumed Density Filtering; Expectation Propagation; 
Kullback-Leibler; Mixture distributions; Probabilistic Editor; Recursive estimation; 
Variational Bayes.
\end{abstract}

\section{Introduction}

Bayesian analysis is typically straightforward if the model from which
the data are generated is compatible with a conjugate family of prior
distributions; this scenario obtains if the model belongs to an
exponential family. However, if there are latent or missing variables
this ideal scenario is not available, even though the model corresponding to the
case in which the latent variables are actually observed 
may be amenable in the sense described
above. In such cases exact calculation of the relevant posterior distributions
and predictive distributions is not available and some form of approximation
becomes inevitable.

Suppose that we have data of the form $x^{(n)}:=
\{x_1,\ldots,x_n\}$, where $n$ denotes the sample size, that the
parameters in the model are denoted by $\theta$ and the corresponding
latent variables are denoted by $z^{(n)}:=
\{z_1,\ldots,z_n\}$. The posterior distribution of interest
is formally given by
$$p(\theta|x^{(n)}) \propto p(x^{(n)}|\theta)~p(\theta),$$
where $p(\theta)$ on the right-hand side is the prior density
for $\theta$. In this paper we shall assume that the data are independently
distributed, so that we can write
$$p(\theta|x^{(n)}) \propto t_0(\theta)\prod^n_{i=1}t_i(\theta),$$
in which $t_0(\theta) = p(\theta)$ and $t_i(\theta) = f(x_i|\theta),$
the probability density or mass function for the $i$th observation, for
$i = 1,\ldots,n$.

One `approximate' approach is to simulate, typically
by some Markov chain Monte Carlo technique,
a large number of realizations from the joint conditional
distribution $p(\theta,z^{(n)}|x^{(n)})$. The resulting realizations
of $\theta$ can be regarded as a sample from $p(\theta|x^{(n)})$
and therefore as a source of empirical approximations of features of interest of
 the true $p(\theta|x^{(n)})$. Similarly, the realizations of the latent
variables provide an approximation to $p(z^{(n)}|x^{(n)})$. In
principle, if enough realizations of $(\theta,z^{(n)})$ are
generated and if the Markov chains that created them can be
guaranteed to have converged, the empirical approximations should
be arbitrarily close to the desired distributions.

However, concerns about the convergence issue and the complexity
of some latent-variable scenarios have led to the development
of deterministic approximations in nontrivial Bayesian contexts.
One approach is to derive a so-called variational approximation, $q$, chosen
to be as close as possible to the true target distribution, $p$,
in terms of the Kullback-Leibler directed divergence
$$\mbox{\sc kl}(q,p) = \int q~{\rm log}~ (q/p),$$
where the `integral' is over $\theta$ and $z^{(n)}$ and has a summation 
component if, as is often the case, the latent variables are discrete.
To make the above optimization practicable some constraints have
to be imposed on $q$, typically that it factorizes into a product of
a function $q_{\theta}$ of $\theta$ and a function $q_z$ of $z^{(n)}$. It turns out
that, if, were the latent variables known, the complete-data model
came from an exponential family and a conjugate prior were used, then
$q_{\theta}$, which represents the approximation
to the intractable posterior distribution of $\theta$, would also belong to 
the conjugate family; see for example Beal and Ghahramani (2003). 
An EM-type algorithm is usually employed to
calculate the appropriate hyperparameters that complete the specification
of the approximating distribution.

While this approach leads to practicable deterministic approximations,
the fact remains that the true posterior distribution is not a member
of the conjugate family; rather, it is a complicated mixture of such
distributions, so that the gap between the true and approximate functions
cannot be made arbitrarily small, for finite sample size, unlike
the case for the simulation-based method. However, there may be some
hope that the variational method will behave acceptably asymptotically,
that is, as the sample size $n$ tends to infinity. In many cases,
in line with the asymptotic properties of maximum likelihood estimators,
posterior distributions tend to be Gaussian, with means ultimately
the same as the maximum likelihood estimates and with covariance
matrices defined in terms of Fisher information matrices. The posterior
mean defines the location of the parameters in the model and a `reliable'
covariance matrix is crucial if interval estimates of the parameters
are required. Thus, any method of approximation can be said to be asymptotically
acceptable if the corresponding $q_{\theta}$ tends to be Gaussian with the correct
mean and correct covariance matrix. A number of papers concerning
mixture problems have established that all is well so far as the
Gaussianity and the mean are concerned, but not so for the covariance
matrix. The limiting covariance matrix tends to that which corresponds
to a complete-data scenario, so that any interval estimates calculated
on the basis of the variational approximation are unrealistically
narrow; see for example Wang and Titterington (2005a, 2005b, 2006).

A different class of deterministic approximations is provided by
the method of Expectation Propagation (Minka, 2001a, 2001b). It is assumed that the
posterior distribution for $\theta$ is approximated by a product
of terms
$$q_{\theta}(\theta) = \prod^n_{i=0}{\tilde t}_i(\theta),$$
where the $i=0$ term corresponds to the prior density and,
typically, if the prior comes from a conjugate family then,
as functions of $\theta$, all the other terms take the same
form so that $q_{\theta}(\theta)$ does also. The explicit form of the
approximation is calculated iteratively. 

{\em Step} 1. From an initial or current
proposal for $q_{\theta}(\theta)$ the $i$th factor ${\tilde t}_i$ is 
discarded and the resulting form is renormalised, giving 
$q_{\theta,{\backslash}i}(\theta)$. 

{\em Step} 2. This $q_{\theta,{\backslash}i}(\theta)$ is then combined with the
`correct' $i$th factor $t_i(\theta)$ (implying that the correct
posterior does consist of a product) and a new $q_{\theta}(\theta)$ of the
conjugate form is selected that gives an `optimal' fit with the new 
product.
In Minka (2001b) optimality is determined by Kullback-Leibler
divergence, in that, if
$$p_i(\theta) \propto q_{\theta,{\backslash}i}(\theta)t_i(\theta),$$
then the new $q_{\theta}(\theta)$ minimises
$$\mbox{\sc kl}(p_i,q) = \int p_i~{\rm log}~ (p_i/q).$$ 
However, in some cases, as we shall see, simpler solutions
are obtained if moment-matching is used instead; if the conjugate
family is Gaussian then the two approaches are equivalent.

{\em Step} 3. Finally, from $q_{\theta,{\backslash}i}(\theta)$ and the new 
$q_{\theta}(\theta)$ a new factor ${\tilde t}_i(\theta)$ can be obtained such that
$$q_{\theta}(\theta) \propto q_{\theta,{\backslash}i}(\theta){\tilde t}_i(\theta).$$

This procedure is iterated over choices of $i$ repeatedly until
convergence is attained. 

As with the variational Bayes approach, an approximate posterior
is obtained that is a member of the complete-data conjugate
family. The important question is to what extent the approximation
improves on the variational approximation. Much empirical
evidence suggests that it is indeed better, and the purpose of this
paper is to investigate this issue at a deeper level. The general 
approach will be to see if Expectation Propagation achieves
what variational Bayes achieves in terms of Gaussianity and
asymptotically correct mean but in addition manages to behave
appropriately in terms of asymptotic second posterior moments.

This paper takes a few first steps by discussing some very simple
scenarios and treating them without complete rigour. In some cases
positive results are obtained, but consideration of one example
suggests that the method is not uniformly successful in the above
terms.\\


\section{Overview of the recursive-based approach to be adopted}

As Minka (2001b) explains, Expectation Propagation was motivated by a recursive
version of the approach known as Assumed Density Filtering (ADF) (Maybeck, 1982). 
In ADF an approximation to the posterior density is created
by incorporating the data one by one, carrying forward a conjugate-family
approximation to the true posterior and updating it observation-by-observation 
using the sort of Kullback-Leibler-based strategy described above in Step 2
of the EP method; the same approach was studied by Bernardo and Gir\'on (1988)
and Stephens (1997, Chapter 5). The difference from EP is that the data are run through only
once, in a particular order, and therefore the final result is order-dependent. 
However, recursive procedures such as these often have desirable asymptotic properties,
sometimes going under the name of stochastic approximations. 
There is also a close relationship with what is called the Probabilistic
Editor; see for example Athans et al. (1977), Makov (1983),
Titterington et al. (1985, Chapter 6) and further references in Section 4. Indeed,
our analysis concentrates on the recursive step defined by Step 2
in the EP approach.

The paper considers in detail normal mixtures with an unknown mean parameter,
for which the conjugate prior family is Gaussian, and mixtures with an unknown
mixing weight, for which the conjugate prior family is Beta. Extension to more
complicated or more general scenarios is currently under investigation.

In this paper we shall deal mainly with one-parameter problems,
although the material in Section 3, at least, can be generalised
easily to the vector case, and we
shall suppose that the conjugate family takes the form $q(\theta|a)$,
where $a$ represents a set of hyperparameters. For simplicity we
shall omit the subscript $\theta$ attached to $q$. Of interest will be
the relationship between consecutive sets of hyperparameters, $a^{(n-1)}$
and $a^{(n)}$, corresponding to the situation before and after the $n$th
observation, $x_n$, is incorporated. Thus, $q(\theta|a^{(n)})$ is the member of
the conjugate family that is closest, in terms of Kullback-Leibler divergence,
or perhaps of moment-matching, to the density that is given by
$$q(\theta|a^{(n-1)})t_n(\theta) = q(\theta|a^{(n-1)})f(x_n|\theta),$$
suitably normalised. If $E_n$ and $V_n$ are the functions of $a^{(n)}$
that represent the mean and the variance of $q(\theta|a^{(n)})$, then
we would want $E_n$ and $V_n$ asymptotically to be indistinguishable from
the corresponding values for the correct posterior.
We would also expect asymptotic Gaussianity.
So far as $E_n$ is concerned, it is a matter of showing that it
converges in some sense to the true value of $\theta$. The correct 
asymptotic variance is essentially given by the asymptotic variance of
the maximum likelihood estimator, with the prior density having
negligible effect, asymptotically. Equivalently, 
the increase in precision associated with the addition of
$x_n$, $V^{-1}_n - V^{-1}_{n-1}$, should be asymptotically the same, in some sense, as
the negative of the second derivative with respect to $\theta$ of 
$\mbox{log}~f(x_n|\theta)$, again by analogy with maximum likelihood theory.

\section{A general finite normal mixture with an unknown mean parameter}

\subsection{Preamble}

As remarked above, only the univariate case is considered although
a multivariate version will follow the same pattern.
This is arguably the most general type of mixture, apart from
regression versions, for which the Gaussian distribution
is the complete-data conjugate prior. The assumption is that the
observed data are a random sample from a univariate mixture of $J$
Gaussian distributions, with means and variances $\{c_j\mu, \sigma^2_j;
j=1,\ldots,J\}$ and with mixing weights $\{v_j; j=1,\ldots,J\}$. The
$\{c_j, \sigma^2_j, v_j\}$ are assumed known, so that $\mu$ is the
only unknown parameter. For observation $x$, therefore, we have
$$f(x|\mu) \propto \sum^J_{j=1}\frac{v_j}{\sigma_j}
\ex\{-\frac{1}{2\sigma^2_j}(x-c_j\mu)^2\}.$$

In the case of a Gaussian distribution the obvious hyperparameters are
the mean and the variance themselves. For comparative simplicity of notation
we shall write those hyperparameters before treatment of the $n$th observation
as $(a, b)$, the hyperparameters afterwards as $(A, B)$ and the $n$th
observation itself as $x$.

\subsection{The recursive step}

The hyperparameters $A$ and $B$ are chosen to match the moments of the
density for $\mu$ that is proportional to
$$b^{-1/2}{\ex}\{-\frac{1}{2b}(\mu-a)^2\}\sum^J_{j=1}\frac{v_j}{\sigma_j}
\ex\{-\frac{1}{2\sigma^2_j}(x-c_j\mu)^2\}.$$
Detailed calculation, described in Appendix 1, shows that the changes in mean
and precision satisfy, respectively,
\begin{equation}
A - a  = b\sum_jR_jS_jT_j/\sum_{j'}{T_{j'}} + o(b)\},
\label{meanupdate}
\end{equation}
and
\begin{equation}
B^{-1}-b^{-1} = \frac{\sum_jR^2_jT_j}{\sum_j{T_{j}}} -
\frac{\sum_jT_jR_j^2S^2_j}{\sum_jT_j} + \frac{(\sum_jT_jR_jS_j)^2}
{(\sum_jT_j)^2} + o(1),
\label{precchange3}
\end{equation}
where $R_j = c_j/\sigma_j$, $S_j = (x-c_ja)/\sigma_j$ and
$$T_j = \frac{v_j}{\sigma_j}\ex\left(-\frac{(x-ac_j)^2}{2\sigma^2_j}\right).$$

\subsection{Fisher information}

As explained earlier, 
the `correct' asymptotic variance is given by the inverse of the
Fisher information, so that the expected change in the inverse
of the variance is the Fisher information corresponding to
one observation, i.e. the negative of the expected second
derivative of
\begin{eqnarray*}
{\rm log}~ f(x|\mu) &=& {\rm const.}+{\rm log}~\left\{\sum^J_{j=1}\frac{v_j}{\sigma_j}
\ex\left(-\frac{1}{2\sigma^2_j}(x-c_j\mu)^2\right)\right\}\\
 & = &  {\rm const.}+{\rm log}~\sum_jT'_j,
\end{eqnarray*}
where $T'_j$ is the same as $T_j$ except that $a$ is replaced by
$\mu$. In what follows, $R'_j$ and $S'_j$ are similarly related
to $R_j$ and $S_j$.
Then it is straightforward to show that
$$\frac{\partial}{{\partial}\mu}~{\rm log}~ f(x|\mu) = \frac{\sum_jT'_jR'_jS'_j}
{\sum_jT'_j}$$
and so the observed information is
\begin{equation}
-\frac{{\partial}^2}{{\partial}\mu^2}~{\rm log}~ f(x|\mu) =
\frac{\sum_jR'^2_jT'_j}{\sum_j{T'_{j}}} -
\frac{\sum_jT'_jR'^2_jS'^2_j}{\sum_jT'_j} + \frac{(\sum_jT'_jR'_jS'_j)^2}
{(\sum_jT'_j)^2}.
\label{info3}
\end{equation}

The right-hand sides of (\ref{precchange3}) and (\ref{info3})
differ only in that (\ref{precchange3}) involves $a$ whereas
(\ref{info3}) involves $\mu$, and (\ref{precchange3}) is
correct just to $O(b)$ whereas (\ref{info3}) is exact.
However, because of the nature of (\ref{meanupdate}),
stochastic approximation theory as applied by 
Smith and Makov (1981) will confirm that asymptotically 
$a$ will converge to $\mu$
and terms of $(O)b$ will be negligible.

Thus, the approximate posterior distribution derived in this section
behaves as we would wish, in terms of its mean and variance;
by construction it is also Gaussian.

\subsection{Precision afforded by the complete-data scenario}

The $\mu$-dependent part of the $n$th observation's contribution to the complete-data
loglikelihood is
$$-\sum_jz_{nj}(x-c_j\mu)^2/(2\sigma^2_j),$$
in which the indicator variable $z_{nj}$ is 1 if the observation
belongs to the $j$th component and is 0 otherwise. The negative
second derivative with respect to $\mu$ is then
$$\sum_jz_{nj}c^2_j/\sigma^2_j = \sum_jR^2_jz_{nj},$$
and, given $x$, $E(z_{nj})=T_j/\sum_{j'}T_{j'}$ so that
\begin{equation}
E(\sum_jz_{nj}c^2_j/\sigma^2_j) = \frac{\sum_jR^2_jT_j}{\sum_j{T_{j}}}.
\label{complete}
\end{equation}

\subsection{Variational Bayes approximation}

It is easy to show that the corresponding change in
precision associated with the variational approximation would give
\begin{eqnarray}
B^{-1}-b^{-1}& = &\frac{\sum_jR^2_jT_j}{\sum_j{T_{j}}} + o(1) \nonumber \\
 & \ge &  \frac{\sum_jR^2_jT_j}{\sum_j{T_{j}}} -
\frac{\sum_jT_jR_j^2S^2_j}{\sum_jT_j} + \frac{(\sum_jT_jR_jS_j)^2}
{(\sum_jT_j)^2} + o(1),
\label{chprec}
\end{eqnarray}
the inequality following because the extra terms can be interpreted as a negative variance.
This indicates that the precision reflected in the variational approximation 
is too high and, as is often the case, is the same as would be achieved in
the complete-data scenario; see (\ref{complete}).

In the next two subsections we consider two particular cases of the
formulation studied so far.

\subsection{A symmetric mixture of two Gaussians}

In this case the model is assumed to correspond to an equally weighted
mixture of $N(-\mu, 1)$ and $N(\mu, 1)$ distributions, where $\mu$
is unknown, so that
$$f(x|\mu) \propto \ex\{-\frac{1}{2}(x+\mu)^2\}+\ex\{-\frac{1}{2}(x-\mu)^2\}.$$
The prior density for $\mu$ is assumed to
be $N(a^{(0)},b^{(0)})$ and the Gaussian approximation for $\mu$ given
data $x^{(n)}:=\{x_1,\ldots,x_n\}$ is assumed to be the $N(a^{(n)},b^{(n)})$
distribution. Consider the incorporation of observation $n$ and, for
simplicity of notation, write $a^{(n-1)}$, $b^{(n-1)}$, $a^{(n)}$,
$b^{(n)}$ and $x_{n}$ as $a$, $b$, $A$, $B$ and $x$, respectively.
This is a particular case of the model studied so far,
with $J=2$, $c_2 = -c_1 = 1$, $\sigma_1 = \sigma_2 = 1$ and $v_1 = v_2 = 1/2.$
Direct substitution gives
\begin{eqnarray*}
 \frac{\sum_jR^2_jT_j}{\sum_j{T_{j}}} &=& 1\\
\frac{\sum_jT_jR_j^2S^2_j}{\sum_jT_j} &=& \frac{(x+a)^2\ex\{-\frac{1}{2}(x+a)^2\}
(x-a)^2\ex\{-\frac{1}{2}(x-a)^2\}}{\ex\{-\frac{1}{2}(x+a)^2\}+\ex\{-\frac{1}{2}(x-a)^2\}}\\
\frac{(\sum_jT_jR_jS_j)^2}{(\sum_jT_j)^2} &=& \left\{\frac
{(x+a)\ex\{-\frac{1}{2}(x+a)^2\}-(x-a)\ex\{-\frac{1}{2}(x-a)^2\}}
{\ex\{-\frac{1}{2}(x+a)^2\}+\ex\{-\frac{1}{2}(x-a)^2\}}\right\}^2,
\end{eqnarray*}
so that (\ref{chprec}) becomes
\begin{equation}
B^{-1}-b^{-1} = 1 - \frac{4x^2\ex\{-\frac{1}{2}(x+a)^2
-\frac{1}{2}(x-a)^2\}}{[\ex\{-\frac{1}{2}(x+a)^2\}+\ex\{-\frac{1}{2}(x-a)^2\}]^2}
+ o(1).
\label{precchange}
\end{equation}
Likewise, the negative of the second derivative of ${\rm log}~ f(x|\mu)$ is
\begin{equation}
-\frac{{\partial}^2}{{\partial}\mu^2}~{\rm log}~ f(x|\mu)
 = 1-\frac{4x^2\ex\{-\frac{1}{2}(x+\mu)^2
-\frac{1}{2}(x-\mu)^2\}}{[\ex\{-\frac{1}{2}(x+\mu)^2\}+\ex\{-\frac{1}{2}(x-\mu)^2\}]^2}.
\label{info}
\end{equation}
As in the general case, the right-hand sides of equations (\ref{precchange}) and (\ref{info})
differ only in that, where the former has $a$, the current posterior mean
for $\mu$, the latter has $\mu$ itself.

Note that the recursion for the posterior mean of the Gaussian approximation
is
$$A = a+b\{(1-w)x-wx-a\}+o(b),$$
where $w = \ex\{-(x+a)^2/2\}/[\ex\{-(x+a)^2/2\}+\ex\{-(x-a)^2/2\}].$
The above equation is similar to standard stochastic approximations 
and should
converge to the true $\mu$, as will the maximum likelihood estimator.

\subsubsection{Quasi-Bayes}

The new observation is assigned to $N(-a,1)$ with probability $w$
and to $N(a,1)$ with probability $1-w$, so that, for appropriate prior,
$$A=a+(n+1)^{-1}\{(1-w)x-wx-a\}$$
and $B=(n+1)^{-1}$. The behaviour
of $b$ in the same as for the confirmed-data case considered next.

\subsubsection{Confirmed-data case}

In this case
$$A=a+(n+1)^{-1}\{-xI(z_n=1)+xI(z_n=2)-a\}$$
and $B=(n+1)^{-1}.$ Also, the Fisher information per observation
is 1. Note that the (observed) information per observation in
(\ref{info}) is less than 1, but approaches 1 as $\mu \rightarrow \infty,$
i.e. as the mixture components become separated.

\subsection{Minka's clutter problem}

\beqn
a &=& b + c \\
  &=& d + e 
\eeqn 

\subsubsection{Preamble}

Again we just consider the univariate version of the problem,
described in Minka (2001b), so that
the model consists of the mixture of $N(\mu, 1)$ and $N(0, 10)$, with known
mixing weights $1-v$ and $v$ respectively; Minka (2001b) uses $w$ for $v$ but we have
 used $w$ already for something else. Thus
$$f(x|\mu) \propto (1-v)\ex\{-\frac{1}{2}(x-a)^2\}+\frac{v}{{\surd}10}~\ex(-x^2/20).$$
As before, we use the notation
$a$, $b$, $A$, $B$ and $x$.
This is another special case of the general model, with $J=2, c_1=1, c_2=0,
\sigma_1 = 1, \sigma_2={\surd}10, v_1=1-v$ and $v_2=v$. Direct substitution gives
\begin{eqnarray*}
 \frac{\sum_jR^2_jT_j}{\sum_j{T_{j}}} &=& \frac{(1-v)\ex\{-\frac{1}{2}(x-a)^2\}}
{(1-v)\ex\{-\frac{1}{2}(x-a)^2\}+\frac{v}{{\surd}10}~\ex(-x^2/20)}\\
\frac{\sum_jT_jR_j^2S^2_j}{\sum_jT_j} &=& \frac{(1-v)(x-a)^2\ex\{-\frac{1}{2}(x-a)^2\}}
{(1-v)\ex\{-\frac{1}{2}(x-a)^2\}+\frac{v}{{\surd}10}~\ex(-x^2/20)}\\
\frac{(\sum_jT_jR_jS_j)^2}{(\sum_jT_j)^2} &=& \frac{[(1-v)(x-a)\ex\{-\frac{1}{2}(x-a)^2\}]^2}
{[(1-v)\ex\{-\frac{1}{2}(x-a)^2\}+\frac{v}{{\surd}10}~\ex(-x^2/20)]^2},
\end{eqnarray*}
so that (\ref{chprec}) becomes
\beq
\begin{array}{rcl}
  B^{-1} - b^{-1}
  &=&\displaystyle \frac{(1-v)\ex\{-\frac{1}{2}(x-a)^2\}}  
  {(1-v)\ex\{-\frac{1}{2}(x-a)^2\}+\frac{v}{{\surd}10}~\ex(-x^2/20)} \\
  & &\displaystyle \qquad - \left(\frac{v(1-v)(x-a)^2
\ex\{-\frac{1}{2}(x-a)^2-x^2/20\}}
{[(1-v)\ex\{-\frac{1}{2}(x-a)^2\}+\frac{v}{{\surd}10}~\ex(-x^2/20)]^2}\right)
+ o(1).
\end{array} 
\label{precchange2}
\eeq 
The negative of the second derivative of ${\rm log}~f(x|\mu)$ is
\beq
\begin{array}{rcl}
-\frac{{\partial}^2~{\rm log}~f(x|\mu)}{{\partial}\mu^2}
&=&\displaystyle 
\frac{(1-v)\ex\{-\frac{1}{2}(x-\mu)^2\}}
     {(1-v)\ex\{-\frac{1}{2}(x-\mu)^2\}+\frac{v}{{\surd}10}~\ex(-x^2/20)} \\
& &\displaystyle \qquad - 
\left(\frac{\frac{v}{{\surd}10}(1-v)(x-\mu)^2
\ex\{-\frac{1}{2}(x-\mu)^2-x^2/20\}}
{[(1-v)\ex\{-\frac{1}{2}(x-\mu)^2\}+\frac{v}{{\surd}10}~\ex(-x^2/20)]^2}\right).
\end{array}
\label{info2}
\eeq 
Again the right-hand sides of equations (\ref{precchange2}) and (\ref{info2})
differ only in that, where the former has $a$, the current posterior mean
for $\mu$, the latter has $\mu$ itself.
Since $a$ will tend to $\mu$ in the
limit, it follows that the approximate approach based on a Gaussian approximation
to the posterior distribution will behave appropriately asymptotically,
so far as mean and variance are concerned.

Note that the recursion for the posterior mean of the Gaussian approximation
is
$$A = a+bw(x-a)+o(b),$$
which is similar to standard stochastic approximations and should
converge to the true $\mu$, as will the maximum likelihood estimator.

\section{Mixture of two known distributions}\label{sec:key}

\subsection{Preamble}

This is arguably the simplest possible mixture but is one for which
the conjugate family is not given by Gaussian distributions.
It is assumed that data arrive independently and identically distributed
from a mixture distribution with density function
$$f(x|\beta) = {\beta}f_1(x) + (1-\beta)f_2(x),$$
in which $\beta$ is an unknown mixing weight between zero and one and
$f_1$ and $f_2$ are known densities. 
The prior density for $\beta$ is assumed to be that of Be$(a^{(0)},b^{(0)}),$
and the Beta approximation for the posterior based on $x^{(n)}:=
\{x_1,\ldots,x_n\}$ is assumed to be the Be$(a^{(n)},b^{(n)})$ distribution,
for hyperparameters $a^{(n)}$ and $b^{(n)}$. The expectation, second moment and variance of
the Beta approximation are respectively
\begin{eqnarray*}
E_n & = & a^{(n)}/(a^{(n)}+b^{(n)}),\\
S_n & = & \{a^{(n)}(a^{(n)}+1)\}/\{(a^{(n)}+b^{(n)})(a^{(n)}+b^{(n)}+1)\}\\
V_n & = & (a^{(n)}b^{(n)})/\{(a^{(n)}+b^{(n)})^2(a^{(n)}+b^{(n)}+1)\}\\
 & = & E_n(1-E_n)/(L_n+1),\\
\end{eqnarray*}
where $L_n = a^{(n)}+b^{(n)}.$
The limiting behaviour of $E_n,$ $S_n$ and $V_n$ is of key interest. 
We would want $E_n$ to tend to the true $\beta$ in some sense and
$V_n$ to tend to the variance of the correct posterior distribution
of $\beta$. Asymptotic normality of the approximating distribution
is also desired.
If $E_n$ behaves as desired then $E_n(1-E_n)$
tends to $\beta(1-\beta)$ and, for $V_n$, the behaviour of $a^{(n)}+b^{(n)}+1
=L_n+1$, and therefore of $L_n$,
is then crucial.

\subsection{The case of confirmed data}

In this case the component identifiers $\{z_n\}$, i.e. $z_n = 1$
or $z_n = 2$, of the observations are known, the exact
posterior distribution is a Beta, and
\begin{eqnarray*}
a^{(n)} & = & a^{(n-1)} + I(z_n=1),\\
b^{(n)} & = & b^{(n-1)} + I(z_n=2),\\
L_n & = & L_{n-1}+1\\
 & = & L_{0}+n\\
   & \bumpeq & n,
\end{eqnarray*}
for large $n$, where $I(\cdot)$ is an indicator function.

In this case $E_n$ is, to first order, ${\hat \beta}_{\rm CO}$, the 
proportion of the $n$  observations that belong to the first component, 
it therefore does tend to $\beta$, by the Law of Large Numbers, 
and the limiting
version of $V_n$ is
$$V_{\rm CO} = \beta(1-\beta)/n$$
for large $n$. Asymptotic normality of the posterior distribution of
$\beta$ follows from the Bernstein-von Mises mirror result corresponding
to the Central Limit result for ${\hat \beta}_{\rm CO}$.

\subsection{The `correct' results}

The behaviour of the `correct' posterior distribution will be
dictated by the behaviour of the maximum likelihood
estimator ${\hat \beta}_{\rm ML}$ which, for this problem, will be consistent
for $\beta$, will be asymptotically normal provided the true
$\beta$ is not 0 or 1 and will have large-sample variance defined by
the Fisher information: for large $n$, approximately,
\begin{eqnarray*}
E({\hat \beta}_{\rm ML}) & = & \beta,\\
\mbox{var}({\hat \beta}_{\rm ML}) & = & \frac{1}{n\int~\frac{
\{f_1(x)-f_2(x)\}^2}{f(x)}dx}\\
 & = & \frac{1}{n\int~\frac{(f_1-f_2)^2}{{\beta}f_1+(1-\beta)f_2}}\\
 & = & V_{\rm ML},
\end{eqnarray*}
say. Again, these properties can be transferred to the posterior
distribution of $\beta$, by a Bernstein-von Mises argument. This transferrence
will apply as a general rule.

\subsection{The variational approximation}

In this case the variational approximation, which is of course a
Beta distribution, has hyperparameters of the form
\begin{eqnarray*}
a^{(n)} &=& a^{(0)} + \sum^n_{i=1} w_{1i}\\
b^{(n)} &=& b^{(0)} + \sum^n_{i=1} w_{2i},
\end{eqnarray*}
where, for each $i$, $w_{1i}$ and $w_{2i}$ are nonzero and sum to 1;
see for example Humphreys and Titterington (2000).
Thus, as in the previous subsection,
$$L_n = a^{(0)}+b^{(0)}+n \bumpeq n,$$
for large $n$. Informally, this implies that the limiting
version of $V_n$ is
$$V_{\rm VA} = \beta(1-\beta)/n.$$
This is the same as $V_{\rm CO}$ and therefore is `smaller than it should be'
and would lead to unrealistically narrow interval estimates for $\beta$.
For more details see for example Humphreys and Titterington (2000)
and, for a more rigorous discussion especially about the convergence
of $E_n$ to $\beta$, Wang and Titterington (2005a).

\subsection{The Quasi-Bayes recursive approach}

For each $n$ let
\begin{equation}
w_{1n} = a^{(n-1)}f_{1n}/(a^{(n-1)}f_{1n}+b^{(n-1)}f_{2n}),
\label{rn}
\end{equation}
where, for $j=1,2,$ $f_{jn} = f_j(x_n).$ Then, for the Quasi-Bayes
approach (Smith and Makov, 1978),
which creates a sequence of Beta approximations that recursively
tracks the expectation from stage to stage,
\begin{equation}
\label{QB}
\begin{aligned}
a^{(n)} & = & a^{(n-1)} + w_{1n},\\
b^{(n)} & = & b^{(n-1)} + 1-w_{1n},
\end{aligned}
\end{equation}
which implies that 
\begin{eqnarray*}
L_n & = & a^{(n-1)}+b^{(n-1)}+1\\
& = & a^{(0)}+b^{(0)}+n\\
& \bumpeq & n.
\end{eqnarray*}
The recursively calculated sequence of posterior means
$\{E_n\}$ satisfy
\begin{equation}
E_n = E_{n-1} + \frac{1}{a^{(n-1)}+b^{(n-1)}+1}(w_{1n} - E_{n-1}),
\label{En}
\end{equation}
Smith and Makov (1978) show, by carefully establishing
the credentials of (\ref{En}) as a stochastic approximation, that the posterior
mean is consistent, so that, for Quasi-Bayes and for large $n$,
$$V_n \bumpeq V_{\rm QB} = \beta(1-\beta)/n;$$
this is the same as for the confirmed-data case and the variational
approximation, thereby falling foul of the same criticism as the latter
as being `too small'.
Similar remarks
apply to any other recursive method in which the hyperparameters
are updated according a rule like that in (\ref{QB}); 
a recursive version of the variational approximation, implemented in
Humphreys and Titterington (2000) is a
case in point.

\subsection{The probabilistic editor}

In Quasi-Bayes the first moment of the posterior distribution
is tracked recursively. In the probabilistic editor (PE) the first
two moments are tracked, which is possible because there are
two hyperparameters. Makov (1983) looks at a context 
different from this simple mixture problem and investigates 
empirically updating more than one observation at a time, but for the
time being we consider just the above framework.

We describe how the hyperparameters are updated based on incorporation
of an observation $x_n$ from the mixture
distribution and a Be$(a^{(n-1)},b^{(n-1)})$ prior distribution. 
For simplicity of notation and in the spirit of previous sections,
we shall denote the hyperparameters by $(a^{(n-1)},b^{(n-1)})=(a,b)$
and $(a^{(n)},b^{(n)})=(A,B)$. Then we have to match moments of
a Be$(A,B)$ distribution with those of the Beta mixture distribution
$$w_1{\rm Be}(a+1,b)+w_2{\rm Be}(a,b+1),$$
where
$$w_1 = af_{1n}/(af_{1n}+bf_{2n}) = 1-w_2.$$
In general, if $\beta \sim \sum_jw_jg_j$, where the distribution
corresponding to $g_j$ has mean $\mu_j$ and variance $\sigma^2_j$, then
\begin{eqnarray*}
E_\beta & = & {\bar \mu} = \sum_jw_j\mu_j,\\
{\rm var}~\beta & =& \sum_jw_j\sigma^2_j + \sum_jw_j(\mu_j
-{\bar \mu})^2.
\end{eqnarray*}
These formulae are also employed in the Appendix.
Thus, matching means, we obtain
$$E_n = \{w_1(a+1)+w_2a\}/(a+b+1)=(a+w_1)/(L_{n-1}+1).$$
Matching variances gives, after some algebra,
$$\frac{E_n(1-E_n)}{L_n+1} = \frac{E_n(1-E_n)}{L_{n-1}+2}
+\frac{w_1(1-w_1)}{(L_{n-1}+2)(L_{n-1}+1)}.$$ 
Clearly, this provides an explicit formula for $L_n$ and consequently for
$A = E_nL_n$ and $B = (1-E_n)L_n,$ but, from a theoretical
point of view, we are interested in approximations that lead to
asymptotic results. Further manipulation gives
$$L_n+1 \bumpeq (L_{n-1}+2)\left\{1-\frac{1}{L_{n-1}+1}\frac{w_1(1-w_1)}
{E_n(1-E_n)}\right\},$$
which then leads to
$$L_n-L_{n-1} \bumpeq 1 - \frac{w_1(1-w_1)}{E_n(1-E_n)} = 1 - \epsilon_n,$$
say, the approximations relying on $L_n$ being of order $n$.
The posterior mean, $E_n$, can be calculated recursively and
 the relevant equation
is the same as (\ref{En}):
\begin{equation*}
E_n = E_{n-1} + \frac{1}{a+b+1}(w_1 - E_{n-1}),
\end{equation*}
where $w_1$ is given in (\ref{rn}). By the argument in Smith and Makov (1978),
the sequence $\{E_n\}$ will converge (to the true $\beta$
value). 

The posterior 
variance of the PE-based Beta approximation based on $n$ observations is
$$V_n = E_n(1-E_n)/(L_n+1),$$
in which the limiting value of $E_n(1-E_n)$ is $\beta(1-\beta).$
For large $n$, $L_n,$ or equivalently
$L_n+1,$ will tend to $nt$, where
$t$ is the limiting expectation of $L_n-L_{n-1}.$ However, for large $n$,
approximately, 
$$w_1(1-w_1) = \beta(1-\beta)f_{1n}f_{2n}/\{{\beta}f_{1n}+(1-\beta)f_{n2}\}^2,$$
and therefore
$$\frac{w_1(1-w_1)}{E_n(1-E_n)} \bumpeq f_{1n}f_{2n}/\{{\beta}f_{1n}+(1-\beta)f_{n2}\}^2.$$
Thus, approximately,
$$E(L_n-L_{n-1}) = 1-\int\frac{f_1f_2}{{\beta}f_1+
(1-\beta)f_2},$$
so that, asymptotically, the variance of the PE-based approximation is
\begin{eqnarray*}
V_{\rm PE} & = & n^{-1}\beta(1-\beta)\left\{1-\int\frac{f_1f_2}{{\beta}f_1+
(1-\beta)f_2}\right\}^{-1}.
\end{eqnarray*}
We now show that $V_{\rm ML} = V_{\rm PE}.$ This follows immediately 
from the following Lemma.

L{\sc emma}. 
\begin{equation}
\frac{1}{\beta(1 - \beta)}\left\{1-\int\frac{f_1f_2}{{\beta}f_1+
(1-\beta)f_2}\right\} = 
\int~\frac{(f_1-f_2)^2}{{\beta}f_1+(1-\beta)f_2}.
\label{peml}
\end{equation}

{\em Proof}. Denote the left-hand and right-hand sides of (\ref{peml})
by $I_1$ and $I_2$, respectively, and note that
$$\beta(1-\beta) = \frac{1}{2}\{1-\beta^2-(1-\beta)^2\}.$$
Then
\begin{eqnarray*}
I_2 &=& \frac{1}{\beta(1 - \beta)}\int\frac{(f^2_1-2f_1f_2+f^2_2)\beta(1-\beta)}
{{\beta}f_1+(1-\beta)f_2}\\
 &=& \frac{1}{\beta(1 - \beta)}\int\frac{\beta(1 - \beta)(f^2_1+f^2_2) + f_1f_2\{\beta^2
+(1-\beta)^2\}-f_1f_2}{{\beta}f_1+(1-\beta)f_2}\\
 &=& \frac{1}{\beta(1 - \beta)}\int\frac{(1-\beta)f_1\{{\beta}f_1+(1-\beta)f_2\}+{\beta}f_2
\{{\beta}f_1+(1-\beta)f_2\}-f_1f_2}{{\beta}f_1+(1-\beta)f_2}\\
&=&\frac{1}{\beta(1 - \beta)}\int\left\{(1-\beta)f_1+{\beta}f_2
-\frac{f_1f_2}{{\beta}f_1+(1-\beta)f_2}\right\}\\
 &=&\frac{1}{\beta(1 - \beta)}\left\{1-\beta+\beta-\int\frac{f_1f_2}{{\beta}f_1+
(1-\beta)f_2}\right\} = I_1.
\end{eqnarray*}

Thus, asymptotically, the Probabilistic Editor, and by implication the
moment-matching version of Expectation Propagation, get the variance right.

\subsection{Results for the Kullback-Leibler update}

The `standard' EP and ADF approaches base the updates not specifically on moment-matching
but on minimization of the Kullback-Leibler divergence discussed in Section 1.
As indicated in Section 3.3.1 of Minka (2001a) and as implied by Section 5.6.4 of
Stephens (1997), the relationships between successive sets of hyperparameters
$(a,b)$ and $(A,B)$ are
\begin{equation}\label{klupdate}
\begin{aligned}
\Psi(A) - \Psi(A+B) &=& \frac{f_{1n}}{af_{1n}+bf_{2n}} -
\frac{1}{a+b} + \Psi(a) - \Psi(a+b)\\
\Psi(B) - \Psi(A+B) &=& \frac{f_{2n}}{af_{1n}+bf_{2n}} -
\frac{1}{a+b} + \Psi(b) - \Psi(a+b),
\end{aligned}
\end{equation}
where $\Psi(c)$ denotes the digamma function. For large $c$, the dominant term
in $\Psi(c)$ is $\mbox{log}~(c-1/2),$ obtainable through Stirling's approximation
to $(c-1)!$, and the corresponding approximations to equations (\ref{klupdate})
are 
\begin{equation}\label{klapprox1}
\begin{aligned}
\mbox{log}~(A-1/2) - \mbox{log}~(A+B-1/2) & = &
\frac{f_{1n}}{af_{1n}+bf_{2n}} - \frac{1}{a+b} + 
\mbox{log}~(a-1/2) - \mbox{log}~(a+b-1/2)\\
\mbox{log}~(B-1/2) - \mbox{log}~(A+B-1/2) & = &
\frac{f_{2n}}{af_{1n}+bf_{2n}} - \frac{1}{a+b} + 
\mbox{log}~(b-1/2) - \mbox{log}~(a+b-1/2).
\end{aligned}
\end{equation}
If we now define $\delta_a$ by $A=a+\delta_a,$
and similarly for $\delta_b$, then we can write
\begin{eqnarray*}
\mbox{log}~(A-1/2)& =& \mbox{log}~(a-1/2)+ \mbox{log}~\{1+\delta_a/(a-1/2)\}\\
 & = & \mbox{log}~(a-1/2) + \delta_a/a + o(a^{-1}),
\end{eqnarray*}
along with similar expressions for $\mbox{log}~(B-1/2)$ and
$\mbox{log}~(A+B-1/2)$. Substitution in (\ref{klapprox1}) gives, for the dominant terms,
\begin{equation}\label{klapprox2}
\begin{aligned}
\delta_a/a-(\delta_a+\delta_b)/(a+b)&=&\frac{f_{1n}}{af_{1n}+bf_{2n}} - \frac{1}{a+b}  \\
\delta_b/b-(\delta_a+\delta_b)/(a+b)&=&\frac{f_{2n}}{af_{1n}+bf_{2n}} - \frac{1}{a+b}.
\end{aligned}
\end{equation}
Substitution of the formulae for $\delta_a$ and $\delta_b$ implicit in
the moment-matching update in Section 4.6 and concentration on
the lowest-order terms lead to equations (\ref{klapprox2})
being satisfied. For example, since $A = L_nE_n = L_n(a+w_{1n})/
(L_{n-1}+1),$ we have
\begin{eqnarray*}
\delta_a/a &=& \{(L_{n}-L_{n-1})a+L_nw_{1n}-a\}/\{a(L_{n-1}+1)\}\\
 & \bumpeq & (L_{n}-L_{n-1})/L_{n-1} + w_{1n}/(aL_{n-1})-1/L_{n-1}\\
 &=& (\delta_a+\delta_b)/(a+b)+f_{1n}/(af_{1n}+bf_{2n})-1/(a+b).
\end{eqnarray*}
The second equation in (\ref{klapprox2}) follows similarly, although
in fact the two equations are linearly dependent.
Thus, to this degree of approximation, the Kullback-Leibler
update is equivalent to the moment-matching update and therefore
performs acceptably, asymptotically.

\subsection{A remark about the above calculation and the treatment of earlier examples}

In the previous examples, in which the conjugate family of distributions
was Gaussian, we found that the (approximate)
increase in precision was `matched' with the single-observation
negative second derivative of the log-density; expectations were
not carried out. The latter seems to be necessary for this example with
an unknown mixing weight, as is now heuristically explained. Let
${\hat \beta}$ denote the estimate of $\beta$ given by the posterior
mean after stage $n-1$, i.e. what has so far been called $E_{n-1}$.
Then the approximate change in precision is
\begin{equation*}
\{{\hat \beta}(1-{\hat \beta})\}^{-1}(L_n-L_{n-1})  \bumpeq 
 \{{\hat \beta}(1-{\hat \beta})\}^{-1}\left(1-\frac{f_{1n}f_{2n}}
{\{{\hat \beta}f_{1n}+(1-{\hat \beta})f_{2n}\}^2}\right),
\end{equation*}
whereas the negative second derivative of the log density for 
the $n$th observation is
\begin{eqnarray*}
-\frac{{\partial}^2}{{\partial}\mu^2}~{\rm log}~ f(x_n|\beta) &=& 
\frac{(f_{1n}-f_{2n})^2}{\{{\beta}f_{1n}+(1-{\beta})f_{2n}\}^2}\\
 &=&\frac{1}{\beta(1-\beta)}\left(\frac
{(1-\beta)f_{1n}+{\beta}f_{2n}}{{\beta}f_{1n}+(1-{\beta})f_{2n}}
-\frac{f_{1n}f_{2n}}
{\{{\beta}f_{1n}+(1-{\beta})f_{2n}\}^2}\right),
\end{eqnarray*}
which does not match up without averaging.

\subsection{Discussion}

In a very small
simulation exercise, Humphreys and Titterington (2000) compare
the non-recursive variational approximation, its recursive variant,
the recursive Quasi-Bayes and Probabilistic Editor, and the Gibbs sampler,
the last of which can be regarded as providing a reliable estimate of the true
posterior. As can be expected on the basis of the above analysis,
the approximation provided by the Probabilistic Editor is very similar
to that obtained from the Gibbs sampler, whereas the other approximations are
`too narrow'. Furthermore, the variances associated with the various
approximations are numerically very close to the `asymptotic' values
derived above.
The implication is that this is also
the case for the corresponding version of the EP algorithm
based on moment matching, mentioned in Section 3.3.3 of Minka (2001a), of which
the PE represents an online version. Of course EP updates using KL divergence
rather than (always) matching moments, but the two versions perform
very similarly in Minka's empirical experiments, and Section
4.7 reflects this asymptotically. Recursive versions of
the algorithm with KL update, i.e. versions of ADF, 
are outlined and illustrated by simulation
in Chapter 5 of Stephens (1997) for mixtures of known densities, 
extending earlier work by Bernardo and Gir{\'o}n (1988),
and for mixtures of Gaussian densities with all parameters unknown,
including the mixing weights. For mixtures of two known densities, Stephens
notes that, empirically, the KL update appears to produce an estimate of
the posterior density that is indistinguishable from the MCMC estimate, and is
much superior to the Quasi-Bayes estimate, which is too narrow. For a mixture
of four known densities, for which the conjugate prior distributions
are 4-cell Dirichlet distributions, the KL update appears clearly to be better than the 
Quasi-Bayes update, but somewhat more `peaked' than it should be. This is because,
in terms of the approach of the present paper, there are insufficient hyperparameters
in the conjugate family to match all first and second moments. For a $J$-cell 
Dirichlet, with $J$ hyperparameters, there are $J-1$ independent first moments 
and $J(J-1)/2$ second moments, so that full moment-matching is not possible for
$J > 2,$ that is, for any case but mixtures of $J=2$ known densities. 
To see this, we proceed along the lines followed in Section 4.6.
Suppose we consider updating a Dir$(a_1, \ldots ,a_J) =$ 
Dir$(a)$ distribution to Dir$(A_1, \ldots ,A_J) =$ Dir$(A)$ 
on the basis of a new observation, $x_n$. Then the
moments of the new Dirichlet have to match those of the Dirichlet mixture
$$\sum_jw_j{\rm Dir}(a+\delta_j),$$
where $a+\delta_j$ denotes the set of hyperparameters that are the same
as the set in $a$ except that the $j$th of them is $a_j+1$ and where, 
for each $j$,
$$w_j = a_jf_{jn}/(\sum_ka_kf_{kn}).$$
Then matching means requires that, for each $j$,
$$E_{jn} = (a_j+w_j)/(L_{n-1}+1) = E_{j,n-1}+\frac{1}{L_{n-1}+1}
(w_{j}-E_{j,n-1}),$$
where $L_{n-1}= \sum_ja_j$, and then
$$A_j = L_{n}E_{jn},$$
where $L_n = \sum_jA_j.$ Clearly, these amount to $J-1$ independent
equations, leaving just one degree of freedom which we assign to the
choice of $L_n.$

Matching variances leads, for each $j$, to 
\begin{equation}
\frac{E_{jn}(1-E_{jn})}{L_n+1} = \frac{E_{jn}(1-E_{jn})}{L_{n-1}+2}
+\frac{w_j(1-w_j)}{(L_{n-1}+2)(L_{n-1}+1)},
\label{matchvar}
\end{equation} 
and then to
\begin{equation}
L_n-L_{n-1} \bumpeq 1 - \frac{w_j(1-w_j)}{E_{jn}(1-E_{jn})}.
\label{approxmv}
\end{equation}
Similarly, matching covariances gives, for each $j$ and $k$ with
$j \ne k,$
\begin{equation}
-\frac{E_{jn}E_{kn}}{L_n+1} = -\frac{E_{jn}E_{kn}}{L_{n-1}+2}
-\frac{w_jw_k}{(L_{n-1}+2)(L_{n-1}+1)},
\label{matchcov}
\end{equation}
and then to
\begin{equation}
L_n-L_{n-1} \bumpeq 1 - \frac{w_jw_k}{E_{jn}E_{kn}}.
\label{approxmc}
\end{equation}
Matching of second-order moments requires $J(J-1)/2$ equations
to be satisfied, and this is clearly possible only for $J=2$.

This is
implicit in the work of Cowell et al. (1996) on recursive updating, 
following on from Spiegelhalter and Lauritzen (1990) and Spiegelhalter 
and Cowell (1992), and referred to in Section 3.3.3 of Minka (2001a)
and Section 9.7.4 of Cowell et al. (1999). 
They chose Dirichlet hyperparameters to match first moments
and the average variance of the parameters. This can clearly be done
exactly, from (\ref{matchvar}), or approximately, from (\ref{approxmv}).
Alternatively one could match the average of the variances and
covariances, based on (\ref{matchcov}) or (\ref{approxmc}).
However, the upshot is that there is no hope that a pure Dirichlet
approximation will produce a totally satisfactory approximation
to the `correct' posterior for $J > 2$, whether through EP or the recursive
alternatives, based on KL updating or moment-matching. However, these versions
should be a distinct improvement on the Quasi-Bayes and variational approximations,
in terms of variance. A possible way forward for small $J$ is to approximate the posterior by a
mixture of a small number of Dirichlets. To match all first- and second-order moments
of a posterior distribution of a set of $J$-cell multinomial probabilities
one would need a mixture of $K$ pure $J$-cell Dirichlets, where
$$KJ+(K-1) = (J-1)+(J-1)+(J-1)(J-2)/2,$$
i.e. Dirichlet hyperparameters + mixing weights = first moments + variances
+ covariances. This gives $K=J/2$. Thus, for even $J$, the match can be exact,
but for odd $J$ there would be some redundancy. In fact, even for $J$ as small
as 4, the algebraic details of the moment-matching become formidable.

In passing we note that there is a directly analogous version of Section
4.7 for the case of $J > 2$.

\bigskip
\noindent{\bf Acknowledgement}

\noindent 
This work has benefited from contact with Philip Dawid, Steffen Lauritzen,
Y.W. Teh and Jinghao Xue. The paper was written, in part while
the authors were in residence at the Isaac Newton Institute
in Cambridge, taking part in the Research Programme on Statistical Theory
and Methods for Complex High-dimensional Data. \\


\section{Appendix: The recursive step for the normal-mixtures problem}

As remarked at the beginning of Section 3.2,
the hyperparameters $A$ and $B$ are chosen to match the moments of the
density for $\mu$ that is proportional to
$$b^{-1/2}{\ex}\{-\frac{1}{2b}(\mu-a)^2\}\sum^J_{j=1}\frac{v_j}{\sigma_j}
\ex\{-\frac{1}{2\sigma^2_j}(x-c_j\mu)^2\}.$$
This turns out to be a Gaussian mixture,
$$\sum_jw_jN(\mu; m_j,s^2_j),$$
where 
\begin{eqnarray*}
w_j & \propto & \frac{v_j}{\sigma_j}\left(\frac{1}{b}+\frac{c^2_j}{\sigma^2_j}
\right)^{-1/2}\ex\left\{-\frac{(x-ac_j)^2}{2b\sigma^2_j\left(
\frac{1}{b}+\frac{c^2_j}{\sigma^2_j}\right)}\right\},\\
m_j &=& \frac{\frac{a}{b}+\frac{c_jx}{\sigma^2_j}}
{\frac{1}{b}+\frac{c^2_j}{\sigma^2_j}},\\
s^2_j &=& \left(\frac{1}{b}+\frac{c^2_j}{\sigma^2_j}\right)^{-1}.
\end{eqnarray*}
Next we match the first two moments with those of a Gaussian
approximation, $N(A,B)$. First, for the mean, we have
\begin{eqnarray*}
A - a &=& \sum_jw_j(m_j-a)\\
 &=& \sum_jw_j\{bR_jS_j + o(b)\},
\end{eqnarray*}
where $R_j = c_j/\sigma_j$ and $S_j = (x-c_ja)/\sigma_j.$
Since in an asymptotic scenario $b$ will be small, so that we shall neglect terms 
of $o(b)$, we only need the first-order, $O(1)$, term in $w_j$:
$$w_j = T_j/\sum_{j'}{T_{j'}},$$
where
$$T_j = \frac{v_j}{\sigma_j}\ex\left(-\frac{(x-ac_j)^2}{2\sigma^2_j}\right),$$
so that 
$$A = a + b\sum_jR_jS_jT_j/\sum_{j'}{T_{j'}} + o(b).$$
Next, for the variance,
$$B = \sum_jw_js^2_j + \sum_jw_jm_j^2 - (\sum_jw_jm_j)^2.$$
Now,
\begin{eqnarray*}
\sum_jw_js^2_j &=& \sum_jw_j\left(\frac{1}{b}+\frac{c^2_j}{\sigma^2_j}\right)^{-1}\\
 & = & \sum_jw_jb\left(1-\frac{bc^2_j}{\sigma^2_j}\right) + o(b^2)\\
 & = & b - b^2\sum_jw_j(c^2_j/\sigma^2_j) + o(b^2).
\end{eqnarray*}
Again it is necessary only to retain the leading term in $w_j$,
so that
$$\sum_jw_js^2_j = b - b^2\sum_jR^2_jT_j/\sum_j{T_{j}} + o(b).$$
Next, if we write
$$m_j = a+\Delta_j,$$
we have that 
\begin{eqnarray*}
\sum_jw_jm^2_j & = & a^2 + 2a\sum_jw_j\Delta_j + \sum_jw_j\Delta^2_j,\\
(\sum_jw_jm_j)^2 & = & a^2 + 2a\sum_jw_j\Delta_j + (\sum_jw_j\Delta_j)^2\\
\end{eqnarray*}
so that, since again only the leading term in $w_j$ is necessary,
\begin{eqnarray*}
\sum_jw_jm_j^2 - (\sum_jw_jm_j)^2 &=& \sum_jw_j\Delta^2_j - (\sum_jw_j\Delta_j)^2\\
 &=& \frac{\sum_jT_j\Delta^2_j}{\sum_jT_j} - \frac{(\sum_jT_j\Delta_j)^2}
{(\sum_jT_j)^2} +o(b^2)\\
 &=& b^2\left(\frac{\sum_jT_jR_j^2S^2_j}{\sum_jT_j} - \frac{(\sum_jT_jR_jS_j)^2}
{(\sum_jT_j)^2}\right) +o(b^2).
\end{eqnarray*}
Thus,
$$B = b\left\{1 - b\left(\frac{\sum_jR^2_jT_j}{\sum_j{T_{j}}} -
\frac{\sum_jT_jR_j^2S^2_j}{\sum_jT_j} + \frac{(\sum_jT_jR_jS_j)^2}
{(\sum_jT_j)^2}\right)\right\} +o(b^2)$$
and therefore
$$B^{-1}-b^{-1} = \frac{\sum_jR^2_jT_j}{\sum_j{T_{j}}} -
\frac{\sum_jT_jR_j^2S^2_j}{\sum_jT_j} + \frac{(\sum_jT_jR_jS_j)^2}
{(\sum_jT_j)^2} + o(1).$$







\medskip

\noindent{\bf References}

\parindent0pt

Athans, M., Whiting, R. and Gruber, M. (1977). A suboptimal estimation algorithm with 
probabilistic editing for false measurements with application to target tracking 
with wake  phenomena. {\em IEEE Trans. Auto. Contr.} {\bf AC-22}, 273--384.

Beal, M.J. and Ghahramani, Z. (2003). The variational Bayesian EM algorithm
for incomplete data: with application to scoring graphical model structures.
In {\em Bayesian Statistics 7}, Ed. J.M. Bernardo, M.J. Bayarri, J.O.
Berger, A.P. Dawid, D. Heckerman, A.F.M. Smith and M. West, pp. 453--464.
Oxford University Press. 

Bernardo, J.M. and Gir\'on, F.J. (1988). A Bayesian analysis of simple
mixture problems (with Discussion). In {\em Bayesian Statistics}, Ed. 
J.M. Bernardo, M.H. DeGroot, D.V. Lindley and A.F.M. Smith, pp. 67--78

Cowell, R.G., Dawid, A.P. and Sebastiani, P. (1996). A comparison of sequential
learning methods for incomplete data. In {\em Bayesian Statistics 5}, 
Ed. J.M. Bernardo, J.O. Berger, A.P. Dawid
and A.F.M. Smith, pp. 533--541. Oxford: Clarendon Press.

Cowell, R.G. Dawid, A.P., Lauritzen, S.L. and Spiegelhalter, D.J. (1999).
{\em Probabilistic Networks and Expert Systems.} New York: Springer.

Humphreys, K. and Titterington, D.M. (2000). Approximate Bayesian inference
for simple mixtures. In {\em COMPSTAT 2000}, Ed. J.G. Bethlehem and P.G.M.
van der Heijden, pp. 331--336. Heidelberg: Physica-Verlag.

Lehmann, E.L. (1991). {\em Theory of Point Estimation.} Pacific Grove, CA:
Wadsworth and Brooks/Cole.

Maybeck, P.S. (1982). {\em Stochastic Models, Estimation and Control.} 
New York: Academic Press.

Makov, U.E. (1983). Approximate Bayesian procedures for dynamic linear models
in the presence of jumps. {\em The Statistician} {\bf 32}, 207--213.

Makov, U.E. and Smith, A.F.M. (1977). A quasi-Bayes unsupervised learning procedure
for priors. {\em IEEE Trans. Inform. Theory} {\bf IT-23}, 761--764.

Minka, T. (2001a). {\em A family of algorithms for approximate Bayesian inference.}
Ph.D. dissertation, Massachusetts Institute of Technology.

Minka, T. (2001b). Expectation Propagation for approximate Bayesian inference.
In {\em Proc. Conf. Uncertainty in AI}.

Owen, R.J. (1975). A Bayesian sequential procedure for quantal response in
the context of adaptive mental testing. {\em J. Am. Statist. Assoc.} {\bf 70},
351--356.

Smith, A.F.M. and Makov, U.E. (1978). A quasi-Bayes sequential procedure for
mixtures. {\em J. R. Statist. Soc.} B {\bf 40}, 106--112.

Smith, A.F.M. and Makov, U.E. (1981). Unsupervised learning for signal
versus noise. {\em IEEE Trans. Inform. Theory} {\bf IT-27}, 498--500.

Spiegelhalter, D.J. and Cowell, R.G. (1992). Learning in probabilistic expert
systems. In {\em Bayesian Statistics 4}, Ed. J.M. Bernardo, J.O. Berger, A.P. Dawid
and A.F.M. Smith, pp. 447--465. Oxford: Clarendon Press.

Spiegelhalter, D.J. and Lauritzen, S.L. (1990). Sequential updating of conditional 
probabilities on directed graphical structures. {\em Networks} {\bf 20}, 579--605.

Stephens, M. (1997). {\em Bayesian Methods for Mixtures of Normal
Distributions.} D. Phil. dissertation, University of Oxford.

Titterington, D.M., Smith, A.F.M. and Makov, U.E. (1985). {\em Statistical Analysis of
Finite Mixture Distributions.} Chichester: Wiley.

Wang, B. and Titterington, D.M. (2005a). Variational Bayes estimation of mixing coefficients.
In {\em Deterministic and Statistical Methods in Machine Learning}, 
Lecture Notes in Artificial Intelligence Vol. 3635, Ed. J. Winkler,
M. Niranjan and N. Lawrence, pp. 281--295. Springer-Verlag.

Wang, B. and Titterington, D.M. (2005b). Inadequacy of interval estimates corresponding to
variational Bayesian approximations. In {\em Proc. 10th Int. Workshop AISTATS}, pp. 373--380.

Wang, B. and Titterington, D.M. (2006). Convergence properties of a general algorithm for 
calculating variational Bayesian estimates for a normal mixture model.
{\em Bayesian Analysis}, {\bf 1} 625--650.

\end{document}